\documentclass[12pt]{article}%
\usepackage{graphicx}
\usepackage{amsmath}
\usepackage{amsfonts}
\usepackage{amssymb}%
\providecommand{\U}[1]{\protect\rule{.1in}{.1in}}
\newtheorem{theorem}{Theorem}
\newtheorem{acknowledgement}[theorem]{Acknowledgement}

\newtheorem{corollary}[theorem]{Corollary}

\newtheorem{definition}[theorem]{Definition}

\newtheorem{proposition}[theorem]{Proposition}
\newtheorem{remark}[theorem]{Remark}

\newenvironment{proof}[1][Proof]{\textbf{#1.} }{\ \rule{0.5em}{0.5em}}
\begin{document}

\title{Multidimensional Chebyshev Systems (Haar systems) - just a definition }
\author{Ognyan Kounchev\\Institute of Mathematics and Informatics\\Bulgarian Academy of Sciences\\and \\IZKS, University of Bonn}
\maketitle

\begin{abstract}
The notion of Chebyshev system was coined by S. Bernstein although many
resutls were proved already by A. Markov. Chebyshev systems play important
role in the one-dimensional Moment problem and Approximation and Spline
theory. We generalize the notion of Chebyshev system for several dimensions
and define the Multidimensional Chebyshev Systems of order $N$. We prove that
this definition is satisfied by the solutions of a wide class of elliptic
equations of order $2N$. This definition generalizes a very large class of
Extended Complete Chebyshev systems in the one-dimensional case. This is the
first of a series of papers in this area, which solves the longstanding
problem of finding a satisfactory multidimensional generalization of the
classical Chebyshev systems introduced already by A. Markov more than hundered
years ago.

\end{abstract}

\begin{center}

\end{center}

\section{History of Chebyshev systems}

\subsection{Developments in the Moment problem by A. Markov and M. Krein}

It was namely in the Moment problem where the notion of Chebyshev systems
appeared for the first time on the big stage, and provided a very natural and
beautiful generalizations of the results of Gauss, Jacobi, Chebyshev,
Stieltjes, Markov, and others.

The classical Moment problem is defined as follows: Find a non-negative
measure $d\mu$ such that
\[
\int_{a}^{b}t^{j}d\mu\left(  t\right)  =c_{j}\qquad\text{for }j=0,1,...,N.
\]
The solution of the problem includes conditions on the constants $c_{j}$
providing solvability. In the case of $N=2n-1$ the problem has been solved by
the famous \emph{Gauss-Jacobi quadrature}; this solution is based on the
orthogonal polynomials $P_{n}$ (of degree $n$ ) which are orthogonal with
respect to the inner product defined by
\[
\left\langle t^{j},t^{k}\right\rangle :=c_{j+k}.
\]
The history is well described in the book of M. Krein and A. Nudelman
\textquotedblright The Markov Moment Problem\textquotedblright%
\ \cite{kreinnudelman}, actually based on a $1951$ paper of M. Krein devoted
to the ideas of Chebyshev. There it is said that A. Markov has realized that
one may consider successfully the Moment problem of the type
\[
\int_{a}^{b}u_{j}\left(  t\right)  d\mu\left(  t\right)  =c_{j}\qquad\text{for
}j=0,1,...,N,
\]
where the system of \textbf{continuous} functions $\left\{  u_{j}\left(
t\right)  \right\}  _{j=0}^{N}$ represent a \textbf{Chebyshev system }in the
interval $\left[  a,b\right]  $, i.e. any linear combination
\[
u\left(  t\right)  =\sum_{j=0}^{N}\alpha_{j}u_{j}\left(  t\right)
\]
has no more than $N$ zeros in $\left[  a,b\right]  $.

Further, by $U_{N}$ we will denote the subspace of the space of continuous
functions $C\left(  \left[  a,b\right]  \right)  $ generated by the Chebyshev
system, i.e.
\begin{equation}
U_{N}:=\left\{  u\left(  t\right)  :u\left(  t\right)  =\sum_{j=0}^{N}%
\alpha_{j}u_{j}\left(  t\right)  \right\}  . \label{UN}%
\end{equation}

In general, in areas other than Approximation theory and Moment problem,
people have tried to find those properties of the one-dimensional algebraic
polynomials which make them so nice. Apparently, the Chebyshev property seems
to be such.

\subsection{Further developments in Approximation theory and Spline theory}

Let us remind also the famous Chebyshev alternance theorem which has been
proved for \textbf{Chebyshev systems}, and which one would like to see in a
multivariate setting (cf. \cite{kreinnudelman}, chapter $9,$ Theorem $4.4$ ) :

\begin{theorem}
(\textbf{Chebyshev-Markov}) Let $f\in C\left(  \left[  a,b\right]  \right)  .$
A necessary and sufficient condition for the element $u_{0}\in U_{N}$ to solve
problem
\[
\inf_{u\in U_{N}}\left\|  f-u\right\|  _{C}=\inf_{u\in U_{N}}\left(
\max_{x\in\left[  a,b\right]  }\left|  f\left(  x\right)  -u\left(  x\right)
\right|  \right)  =:\delta
\]
is the existence of $N+2$ points
\[
t_{1}<\cdot\cdot\cdot<t_{N+2}%
\]
such that
\[
\delta\varepsilon\left(  -1\right)  ^{j}=f\left(  t_{j}\right)  -u_{0}\left(
t_{j}\right)  \qquad\text{for }j=1,...,N+2
\]
where $\varepsilon=1$ or $\varepsilon=-1.$
\end{theorem}

What concerns other areas where Chebyshev systems have found numerous
applications, one has to mention the book of \cite{schumaker} which contains
an exhaustive consideration of spline theory where splines are piecewise
functions belonging to a Chebyshev system.

\section{Definitions}

Let us provide some basic definitions.

Consider the system of functions $\left\{  u_{j}\left(  t\right)  \right\}
_{j=0}^{N}$ defined on some interval $\left[  a,b\right]  $ in $\mathbb{R}$
and the linear space defined in (\ref{UN}).

\begin{definition}
\label{DTsystem}We call the system of functions $\left\{  u_{j}\left(
t\right)  \right\}  _{j=0}^{N}$ Chebyshev (or $T-$system) iff for every set of
constants $\left\{  c_{j}\right\}  _{j=0}^{N}$ and every choice of the points
$t_{j}\in\left[  a,b\right]  $ with
\[
t_{0}<t_{1}<\cdot\cdot\cdot<t_{N}%
\]
there is a unique solution $u\in U$ of the equations
\[
u\left(  t_{j}\right)  =c_{j}\qquad\text{for }j=0,1,...,N.
\]

It is equivalent to say that
\[
u\left(  t_{j}\right)  =0\qquad\text{for }j=0,1,...,N
\]
implies%
\[
u\equiv0.
\]

\end{definition}

\begin{proposition}
Assume that the space $U\subset C\left[  a,b\right]  $ is given. Then if for
some set of knots $t_{j}\in\left[  a,b\right]  $ with
\[
t_{0}<t_{1}<\cdot\cdot\cdot<t_{N},
\]
and for arbitrary constants $\left\{  c_{j}\right\}  $ we have unique solution
$u\in U$ of the equations
\[
u\left(  t_{j}\right)  =c_{j}\qquad\text{for }j=0,1,...,N
\]
it follows that $\dim U=N+1.$
\end{proposition}

One may formulate the above in an equivalent way:

\begin{proposition}
The following are equivalent

1. the system $\left\{  u_{j}\left(  t\right)  \right\}  _{j=0}^{N}$ is $T-$system

2. for every $u\in U_{N}$ the number of zeros in the interval $\left[
a,b\right]  $ is $\leq N.$

3. the following determinants satisfy
\[
D\left(  u;t_{0},t_{1},...,t_{N}\right)  :=\det\left[
\begin{array}
[c]{cccc}%
u_{0}\left(  t_{0}\right)  & u_{0}\left(  t_{1}\right)  & \cdot\cdot\cdot &
u_{0}\left(  t_{N}\right) \\
u_{1}\left(  t_{0}\right)  & u_{1}\left(  t_{1}\right)  & \cdot\cdot\cdot &
u_{1}\left(  t_{N}\right) \\
\cdot & \cdot & \cdot & \cdot\\
u_{N}\left(  t_{0}\right)  & u_{N}\left(  t_{1}\right)  & \cdot\cdot\cdot &
u_{N}\left(  t_{N}\right)
\end{array}
\right]  \neq0
\]

\end{proposition}

\begin{definition}
\label{DT+} If $D\left(  u;t_{0},t_{1},...,t_{N}\right)  \geq0$ then $\left\{
u_{j}\left(  t\right)  \right\}  _{j=0}^{N}$ is called $T_{+}-$system.
\end{definition}

\subsection{Examples}

The classical polynomials, the trigonometric polynomials in smaller intervals
$\left[  0,2\pi\right]  $!

\begin{enumerate}
\item the system
\[
\left\{  u_{j}\left(  t\right)  =t^{\alpha_{j}}\right\}  _{j=0}^{N}%
\qquad\text{on subintervals of }\left[  0,\infty\right]
\]

\item the system
\[
u_{j}\left(  t\right)  =\frac{1}{s_{j}+t}\qquad\text{for }0<s_{0}<s_{1}%
<\cdot\cdot\cdot<s_{N}\qquad\text{on closed subint. of }\left(  0,\infty
\right)  .
\]

\item the system
\[
u_{j}\left(  t\right)  =e^{-\left(  s_{j}-t\right)  ^{2}}\qquad\text{for
}0<s_{0}<s_{1}<\cdot\cdot\cdot<s_{N}\qquad\text{on }\left(  -\infty
,\infty\right)  .
\]

\item if $G\left(  s,t\right)  $ is the Green function associated with the
operator
\[
Lf=-\frac{d}{dx}\left(  p\frac{df}{dx}\right)  +qf
\]
and some boundary conditions on the interval $\left[  a,b\right]  ,$ then the
system
\[
u_{j}\left(  t\right)  =G\left(  s_{j},t\right)  \qquad\text{for }%
0<s_{0}<s_{1}<\cdot\cdot\cdot<s_{N}\qquad\text{on closed subint. of }\left[
a,b\right]
\]

\end{enumerate}

For further examples, see the monograph of M. Krein and A. Nudel'man
\cite{kreinnudelman} and of S. Karlin and W. Studden, \cite{karlinstudden}.

\section{Extended systems}

We usually work with differentiable systems of functions, and we count the
multiplicities of the zeros.

\begin{definition}
Let $\left\{  u_{j}\left(  t\right)  \right\}  _{j=0}^{N}\in C^{N}\left[
a,b\right]  $ be a $T-$system. We call it Extended Chebyshev system
($ET-$system) if in $U_{N}$ we may uniquely solve the \textbf{Hermite
interpolation problem}:%
\[
u^{\left(  k\right)  }\left(  t_{j}\right)  =c_{j,k}\qquad\text{for
}k=0,1,...,d_{j}%
\]
with arbitrary numbers $\left\{  c_{j,k}\right\}  $ where
\[
\sum\left(  d_{j}+1\right)  =N+1.
\]

It is equivalent to say that if for some $u\in U_{N}$ holds
\[
u^{\left(  k\right)  }\left(  t_{j}\right)  =0\qquad\text{for }k=0,1,...,d_{k}%
\]
then
\[
u\equiv0.
\]

\end{definition}

There are equivalent formulations with determinants and zeros:

\begin{proposition}
The following are equivalent:

1. the system $\left\{  u_{j}\left(  t\right)  \right\}  _{j=0}^{N}\in
C^{N}\left[  a,b\right]  $ is and $ET-$system

2. for every $u\in U$ the number of zeros counted with the multiplicities is
$\leq N.$

3. the \textbf{modified} determinants (see the detailed explanation in
\cite{kreinnudelman}, end of section $1,$ chapter $2$, or \cite{karlinstudden}%
),
\[
\det\left[
\begin{array}
[c]{ccccccc}%
u_{0}\left(  t_{0}\right)  & u_{0}^{\prime}\left(  t_{0}\right)  & \cdot
\cdot\cdot & u_{0}^{\left(  k_{1}-1\right)  }\left(  t_{0}\right)  &
u_{0}\left(  t_{1}\right)  & \cdot\cdot\cdot & u_{0}\left(  t_{N}\right) \\
u_{1}\left(  t_{0}\right)  & u_{1}^{\prime}\left(  t_{0}\right)  & \cdot
\cdot\cdot & u_{1}^{\left(  k_{1}-1\right)  }\left(  t_{0}\right)  &
u_{1}\left(  t_{1}\right)  & \cdot\cdot\cdot & u_{1}\left(  t_{N}\right) \\
\cdot & \cdot & \cdot & \cdot & \cdot & \cdot & \cdot\\
u_{N}\left(  t_{0}\right)  & u_{N}^{\prime}\left(  t_{0}\right)  & \cdot
\cdot\cdot & u_{N}^{\left(  k_{1}-1\right)  }\left(  t_{0}\right)  &
u_{N}\left(  t_{1}\right)  & \cdot\cdot\cdot & u_{N}\left(  t_{N}\right)
\end{array}
\right]  \neq0
\]

\end{proposition}

There is a nice \textbf{characterization} of $ET-$systems. A basic example of
$ET-$system is the following: Let $w_{i}\in C^{N-i}\left[  a,b\right]  $ be
positive functions on $\left[  a,b\right]  $ for $i=0,1,...,N.$ Then the
functions
\begin{align}
u_{0}\left(  t\right)   &  =w_{0}\left(  t\right)  \label{canonical1}\\
u_{1}\left(  t\right)   &  =w_{0}\left(  t\right)  \int_{a}^{t}w_{1}\left(
t_{1}\right)  dt_{1}\\
u_{2}\left(  t\right)   &  =w_{0}\left(  t\right)  \int_{a}^{t}w_{1}\left(
t_{1}\right)  \int_{a}^{t_{1}}w_{2}\left(  t_{2}\right)  dt_{2}\cdot dt_{1}\\
&  \cdot\cdot\cdot\\
u_{N}\left(  t\right)   &  =w_{0}\left(  t\right)  \int_{a}^{t}w_{1}\left(
t_{1}\right)  \int_{a}^{t_{1}}w_{2}\left(  t_{2}\right)  \cdot\cdot\cdot
\int_{a}^{t_{N-1}}w_{N}\left(  t_{N}\right)  dt_{N}\cdot\cdot\cdot
dt_{1}\label{canonical2}%
\end{align}
form an $ET-$system.

\begin{definition}
\label{DM+} An $T-$system $\left\{  u_{k}\right\}  _{k=0}^{N}$ is called
$M_{+}-$system if every subsystem $\left\{  u_{k}\right\}  _{k=0}^{m}$ for
$m=0,1,...,N$ is a $T_{+}-$system.
\end{definition}

There is an important result which is basically due to S. Bernstein (cf.
\cite{kreinnudelman}, Theorem $4.1$ and  Theorem $4.2,$ chapter $2,$ and the
footnote at the end).

\begin{theorem}
Every $T-$ system $\left\{  u_{k}\right\}  _{k=0}^{N}$ on the interval
$\left[  a,b\right]  $ may be linearly transformed to an $M_{+}-$system
$\left\{  v_{k}\right\}  _{k=0}^{n}$ on the interval $\left(  a,b\right)  $.
\end{theorem}

A similar result may be proved in the differentiable case. 

\begin{corollary}
Every $ET-$system in $\left[  a,b\right]  $ may be linearly transformed to an
$ECT-$system in $\left(  a,b\right)  .$
\end{corollary}

The role of the $ECT$ systems becomes clear from the following fundamental result.

\begin{theorem}
\label{T4.2and5.1} If the space $U_{N}$ is generated by an $ECT-$system of
order $N$ then it has a basis $\left\{  v_{j}\right\}  _{j=0}^{N}$ which is
representable in the above form (\ref{canonical1})-(\ref{canonical2}). Hence,
$U_{N}$ is a set of solutions to the following equation
\[
L_{N}u=0\qquad\text{for }t\in\left(  a.b\right)
\]
and
\begin{equation}
L_{N}=\prod_{j=0}^{N}\frac{d}{dt}\frac{1}{w_{j}\left(  t\right)  }.\label{LN}%
\end{equation}

\end{theorem}

The proof is available in \cite{karlinstudden} (chapter $11,$ Theorem $1.1$)
and for $M_{+}-$ systems in \cite{kreinnudelman}, Theorems $4.1$,  $4.2$ and
$5.2$, chapter $2$.

\section{The multivariate case -- attempts}

We have to note that all "brute force generalizations" of the Chebyshev
systems fail. We will mention some of them.

\subsection{Generalization by zero sets -- theorem of \textbf{Mairhuber}
\label{SMairhuber}}

Apparently, the first attempt has been to generalize the Chebyshev systems by
considering the set of zeros:

\begin{definition}
\label{Dzeros}Let $K$ be a compact topological space. The system of functions
$\left\{  u_{j}\right\}  _{j=0}^{N}$ is called Chebyshev of order $N$ iff
\[
Z\left(  u\right)  \leq N
\]
for every $u\in U.$
\end{definition}

The following result shows that there are no non-trivial examples of
multidimensional systems satisfying Definition \ref{Dzeros} (cf.
\cite{kreinnudelman}, chapter $2,$ section $1$).

\begin{theorem}
\label{TMairhuber} (Mairhuber, $1956$) The only spaces $K$ having a Chebyshev
system satisfy $K\subset\mathbb{R}$ or $K\subset\mathbb{S}^{1}$ .
\end{theorem}

The result of Theorem \ref{TMairhuber} is intuitively clear since "general
position" function in $C\left(  K\right)  $ has a zero set which is a subset
of $K$ of codimension $1.$ In particular, if $K=\mathbb{R}^{2}$ then the zero
set is roughly speaking union of some curves, and it would be more reasonable
to speculate about the number of these components then to consider Definition
\ref{Dzeros} above. Speculating in this direction in the multidimensional
case, in view of the Polyharmonic Paradigm \cite{okbook}, one may try to
replace the points on $\mathbb{R}$ by closed surfaces,  and the spaces $U$ by
solutions of Elliptic PDEs. Going further, one may obtain some interesting
results if one uses spheres in the case of the polyharmonic operator
$\Delta^{N}$ by considering the space
\[
U_{N}=\left\{  u:\Delta^{N}u=0\quad\text{in }D\right\}
\]
for some bounded domain $D\subset\mathbb{R}^{n}.$  In particular, one may
prove that if a polyharmonic function of order $N$ in a domain $D$ (a function
satisfying $\Delta^{N}u=0$ in $D$) is zero on a set of $n$ concentric spheres,
then $u\equiv0.$ This result has been proved apparently a long time ago by
means of the Almansi theorem, see  e.g. \cite{atakhodzaev}. However these
hopes to try to generalize Definition \ref{Dzeros} are only vain. They  have
been broken by and example which has been published apparently for a first
time in $1982$, by Atakhodzhaev, \cite{atakhodzaev}. It shows a non-zero
biharmonic function which is zero on two embedded ovals in $\mathbb{R}^{2}.$ 

\subsection{Generalization by Haar property}

The following result belongs to A. Haar, \cite{kreinnudelman}:\ 

\begin{theorem}
Let the space $U_{N}\subset C\left[  a,b\right]  $ be generated by a Chebyshev
system $\left\{  u_{j}\right\}  _{j=0}^{N}.$ Then for every $f\in C\left[
a,b\right]  $ the best approximation problem
\[
\min_{u\in U_{N}}\left\Vert f-u\right\Vert _{C}%
\]
has unique solution.
\end{theorem}

Extending this definition to the multivariate case seems to be very reasonable
but the work with best approximations is very heavy and until now has not led
to success.

\section{Systems to be generalized}

One needs a new point of view on the Chebyshev systems which would make them
generalizable to several dimensions. We propose the point of view of boundary
value problems: We consider a special class of $ET-$systems which are ''generalizable''.

\begin{definition}
\label{DDTsystem}We say that the system $\left\{  u_{j}\right\}  _{j=0}%
^{2N-1}\in C^{2N-1}\left[  a,b\right]  $ is a \textbf{Dirichlet type Chebyshev
system,} or $DT-$system, in the interval $\left[  a,b\right]  $ if for every
two points $\alpha$ and $\beta$ in $\left[  a,b\right]  $ and for every set of
constants $c_{j}$ and $d_{j}$ we are able to \textbf{solve uniquely} the
following \textbf{interpolation problem}, with $u\in U_{N},$
\begin{align*}
u^{\left(  k\right)  }\left(  \alpha\right)   &  =c_{k}\qquad\text{for
}k=0,1,...,N-1\\
u^{\left(  k\right)  }\left(  \beta\right)   &  =d_{k}\qquad\text{for
}k=0,1,...,N-1
\end{align*}

\end{definition}

\begin{remark}
Obviously, all $ET-$systems are $DT-$systems but not vice versa.
\end{remark}

Let us state an equivalent formulation which we are going to  mimic in the
multivariate case.

\begin{proposition}
\label{PDirichlet}The system $\left\{  u_{j}\right\}  _{j=0}^{2N-1}\in
C^{\infty}\left[  a,b\right]  $ is a \textbf{Dirichlet type Chebyshev system}
in the interval $\left[  a,b\right]  $ iff for every two points $\alpha$ and
$\beta$ in $\left[  a,b\right]  ,$  for every set of constants $c_{j}$ and
$d_{j},$ and for every $\varepsilon>0,$ we are able to \textbf{solve } the
following approximate interpolation problem, with $u\in U_{2N}%
=\operatorname*{span}\left\{  u_{j}\right\}  _{j=0}^{2N-1},$
\begin{align}
\left\vert u^{\left(  k\right)  }\left(  \alpha\right)  -c_{k}\right\vert  &
<\varepsilon\qquad\text{for }k=0,1,...,N-1\label{DirichletApprox1}\\
\left\vert u^{\left(  k\right)  }\left(  \beta\right)  -d_{k}\right\vert  &
<\varepsilon\qquad\text{for }k=0,1,...,N-1\label{DirichletApprox2}%
\end{align}

\end{proposition}

The proof is evident since the space $U_{2N}$ is finite-dimensional.

\section{The multivariate case}

In the present section we will provide a multivariate generalization to the
$DT$-systems of Definition \ref{DDTsystem}. 

One might use the properties of the Dirichlet type Chebyshev system provided
in Proposition \ref{PDirichlet} as possible way to make a multivariate
generalization. However we would like to have also the properties exposed by
Theorem \ref{T4.2and5.1} as well. 

\begin{remark}
It is expected that a further research would prove equivalence between the
solvability of problem (\ref{DirichletApprox1})-(\ref{DirichletApprox2}) and
the representation of the space $U_{N}$ as a set of solutions to an equation
$L_{N}u=0.$ 
\end{remark}

\textbf{Whatever the definition of \emph{Multidimensional Chebyshev systems},
we would like to retain} the properties in Theorem \ref{T4.2and5.1} and
Proposition \ref{PDirichlet}. In general, it would need in the future to make
a proper refinement of these properties which would make then into two
equivalent sets of conditions. However at the present moment we will restrict
ourselves with some special though sufficiently wide generalization. 

We consider a subspace of functions $U$ with $U\subset C^{\infty}\left(
D\right)  $ for some bounded domain $D\subset\mathbb{R}^{n}$ such that its
boundary $\partial D$ is infinitely smooth,  and assume that $D$ locally
\textquotedblright lies on one side of the boundary\textquotedblright. These
are the usual conditions for the solvability of Elliptic Boundary Value
problems, see e.g. \cite{lionsmagenes}. For simplicity assume that $D$ is
connected and  simply connected as well.

In the following definition  we will mimic the properties of the Chebyshev
systems provided  in Theorem \ref{T4.2and5.1}. 

\begin{definition}
\label{Dfactorizable} We will say that the elliptic operator $P_{2N}$ of order
$2N$ defined in the domain $D,$  is \textbf{factorizable} if there exist $N$
uniformly strongly elliptic operators $Q_{2}^{\left(  j\right)  }$ of second
order,  defined in the domain $D,$ and  satisfying the following properties:

1. Every operator $Q_{2}^{\left(  j\right)  }$ satisfies the maximum principle
in $\overline{D}.$ 

2. Every operator $Q_{2}^{\left(  j\right)  }$ satisfies condition $\left(
U\right)  _{s}$ for uniqueness in the Cauchy problem in the
small.\footnote{The differential operator $P$ satisfies condition $\left(
U\right)  _{s}$ for uniqueness in the Cauchy problem in the small in $G$
provided that if $G_{1}$ is a connected open subset of $G$ and $u\in
C^{r}\left(  G_{1}\right)  $ is a solution to $P^{\ast}u=0$ and $u$ is zero on
a non-emplty subset of $G_{1}$ then $u$ is identically zero. Elliptic
operators with analytic coefficients satisfy this property (cf.
\cite{bersJohnschechter}, part $II,$ chapter $1.4;$ \cite{browder}, p.
$402$).}

3. The following equality holds
\[
P_{2N}=\prod_{j=1}^{N}Q_{2}^{\left(  j\right)  }=Q_{2}^{\left(  1\right)
}Q_{2}^{\left(  2\right)  }\cdot\cdot\cdot Q_{2}^{\left(  N\right)  }.
\]

\end{definition}

However in the next definition we will mimic the interpolation  properties of
the Chebyshev systems exposed by Proposition \ref{PDirichlet}. 

\begin{definition}
\label{DMultivariateDT} We say that the space  $U$ satisfies the
\textbf{Multivariate BVP Interpolation} of order $N$ iff the following
conditions hold:\ 

1. The  \textbf{approximate solvability of BVP}  on subdomains holds in the
following sense: Let the "boundary differential operators" $B_{j}\left(
x;D_{x}\right)  ,$ $j=1,2,...,N,$ with smooth coefficients and of orders
$\leq2N,$ defined in $\overline{D}$ be given. Let $D_{1}$ be an arbitrary
subdomain of $D$ with $D_{1}\subset D,$ and such that $D_{1}$ satisfies the
above conditions as $D$, and also  $D\setminus D_{1}$ has only non-compact
connected components. Let $c_{j}\in C^{\infty}\left(  \partial D_{1}\right)  $
for $j=0,1,...,N-1$. Then for every $\varepsilon>0$ there exists an element
$u\in U$ such that
\begin{equation}
\left\vert B_{j}u-c_{j}\left(  x\right)  \right\vert \leq\varepsilon
\qquad\text{for all }x\in\partial D_{1},\text{and }%
j=1,...,N.\label{ApproximateInterpolation}%
\end{equation}

In the case of the whole domain, i.e.  $D_{1}=D,$  inequality
(\ref{ApproximateInterpolation}) holds with   $\varepsilon=0.$

2. The \textbf{unique solvability of BVP} on subdomains holds: If for some
$u\in U$ holds
\begin{equation}
B_{j}u\left(  x\right)  =0\qquad\text{for all }x\in\partial D_{1},\text{ for
}j=1,...,N,\label{Uniqueness}%
\end{equation}
then
\[
u\equiv0.
\]

\end{definition}

\begin{remark}
In the case of the domain $D$ we have unique solvability of the Elliptic BVP
since we may take $\varepsilon=0$ !
\end{remark}

\begin{remark}
For Theorem \ref{Tmain} below it is important to note that if an operator
$B_{j}$ is non-characteristic in every direction then it is elliptic, see
\cite{lionsmagenes}, (Definition $1.4$ in chapter $2,$ section $1.4$ ). 
\end{remark}

We would like for a \textbf{Multidimensional Chebyshev space} $U$ to satisfy
analogs to both Theorem \ref{T4.2and5.1} and Proposition \ref{PDirichlet}. In
this respect we may prove the following theorem which shows that it is better
to define the \textbf{Multidimensional Chebyshev spaces }by means of
Definition \ref{Dfactorizable} than by means of Definition
\ref{DMultivariateDT}. 

\begin{theorem}
\label{Tmain} Let the elliptic operator $P_{2N}$ in the domain $D$ be
\textbf{factorizable} by Definition \ref{Dfactorizable}. Then the space
\[
U_{N}=\left\{  u\in H^{2N}\left(  D\right)  :P_{2N}u=0\quad\text{in
}D\right\}
\]
satisfies the \textbf{Multivariate BVP interpolation} of Definition
\ref{DMultivariateDT}.
\end{theorem}%

\proof
We choose the following boundary operators
\[
B_{j}=\prod_{i=N-j+2}^{N}Q_{2}^{\left(  i\right)  }\qquad\text{for }j\geq1,
\]
and in particular, $B_{1}=id.$ 

The uniqueness (\ref{Uniqueness}) follows by induction in $N,$ from item $1$
in Definition \ref{Dfactorizable}. 

Item $1$ in Definition \ref{DMultivariateDT} follows inductively in $N.$ Let
us consider for simplicity the case $N=2.$ Let $u_{0}\in H^{4}\left(
D_{1}\right)  $ satisfy $P_{4}u_{0}=Q^{\left(  1\right)  }Q^{\left(  2\right)
}u_{0}=0$ in $D_{1}.$ Let us put
\[
w_{0}=Q^{\left(  2\right)  }u_{0}%
\]
By the $\left(  U\right)  _{s}$ property of the operator $Q^{\left(  1\right)
},$ it follows that a \textbf{Runge} type theorem holds, namely, for every
$\varepsilon>0$ there exists a solution $w_{\varepsilon}\in H^{2}\left(
D\right)  $ of the equation $Q^{\left(  1\right)  }w_{\varepsilon}=0$ in $D$
and
\[
\left\Vert w_{\varepsilon}-w_{0}\right\Vert \leq\varepsilon,
\]
(cf. \textbf{ }\cite{lax}, \cite{malgrange}, \cite{browder}). Now we want to
find a solution $u_{\varepsilon}\in H^{4}\left(  D\right)  $ such that
$Q^{\left(  2\right)  }u_{\varepsilon}=w_{\varepsilon}$ in $D$ and $\left\Vert
u_{\varepsilon}-u_{0}\right\Vert _{D_{1}}\leq2\varepsilon.$ But here we use
the $\left(  U\right)  _{s}$ property of the operator $Q^{\left(  2\right)  }$
since we compare $Q^{\left(  2\right)  }u_{0}=w_{0}$ and $Q^{\left(  2\right)
}u_{\varepsilon}=w_{\varepsilon}$ and we know that $\left\Vert w_{\varepsilon
}-w_{0}\right\Vert \leq\varepsilon.$ By the $\left(  U\right)  _{s}$ property
as above we may find a solution $u_{\varepsilon}$ to the non-homogeneous
equation which satisfies $\left\Vert u_{\varepsilon}-u_{0}\right\Vert _{D_{1}%
}\leq2\varepsilon.$ 

This ends the proof. %

\endproof

\begin{remark}
An alternative reference for the Runge-Lax-Malgrange  type theorem is
\cite{hoermander1} (Theorem $4.4.5$), \cite{tarkhanov} ( Section $2.3,$
Theorems $2.1,$ $2.2,$ and Theorem $2.4$ for which it is mentioned there, that
it was formulated by F. Browder with error). 
\end{remark}

\section{Examples}

For some integer $N\geq1$ let us consider the space
\[
U:=\left\{  u:\Delta^{N}u\left(  x\right)  =0\qquad\text{in }D\right\}  .
\]

\begin{proposition}
The space $U$ satisfies  Definition \ref{Dfactorizable}.  
\end{proposition}

The proof follows directly after we  define the operators
\[
B_{j}=\Delta^{j-1}\qquad\text{for }j=1,2,...,N.
\]

\subsection{The one-dimensional case}

We still have to check that the $DT-$systems defined in Definition
\ref{DDTsystem} satsfy the multivariate counterpart in Definition
\ref{DMultivariateDT}.

\begin{proposition}
Let us assume that in Definition \ref{DMultivariateDT} the space dimension is
$n=1.$ Then the set $U$ coincides with a  \textbf{Dirichlet type Chebyshev
system} $DT$ of order $N$ from Definition \ref{DDTsystem}. 
\end{proposition}%

\proof
Indeed, let us take the set $D=\left[  a,b\right]  $ and apply the
interpolation propety to the case $D_{1}=D.$ Then we know that for all
constants $c_{j}$ and $d_{j}$ we have unique solvability of the problem
\begin{align*}
u^{\left(  k\right)  }\left(  a\right)   &  =c_{k}\qquad\text{for
}k=0,1,...,N-1\\
u^{\left(  k\right)  }\left(  b\right)   &  =d_{k}\qquad\text{for
}k=0,1,...,N-1.
\end{align*}
Hence, $U$ is $2N$ dimensional; we may take the solution $v_{j}$ to the
problem
\begin{align*}
u^{\left(  k\right)  }\left(  a\right)   &  =\delta_{j,k}\qquad\text{for
}k=0,1,...,N-1\\
u^{\left(  k\right)  }\left(  b\right)   &  =0\qquad\text{for }k=0,1,...,N-1
\end{align*}
and the solution $w_{j}$ to the problem
\begin{align*}
u^{\left(  k\right)  }\left(  a\right)   &  =0\qquad\text{for }k=0,1,...,N-1\\
u^{\left(  k\right)  }\left(  b\right)   &  =\delta_{j,k}\qquad\text{for
}k=0,1,...,N-1
\end{align*}
and we make a basis for $U.$ The approximate solvability of the Dirichlet
problem for $D_{1}\subset\left[  a,b\right]  $ implies now the exact
solvability since $U$ is finite-dimensional. Indeed, we will take a sequence
of solutions $u_{\varepsilon}\left(  t\right)  $ and the limit.%

\endproof

\section{Remarks}

\begin{enumerate}
\item In the case of Multidimensional Chebyshev Systems we need the
approximate solvability of the Dirichlet problem since we have
infinite-dimensional spaces, and there is no equivalence between the
uniqueness and the existence but we have a substitute which is the Fredholm
property of the regular  Elliptic BVPs. The one-dimensional Proposition
\ref{PDirichlet} traces the smooth path for the multivariate generalization. 
\end{enumerate}

\begin{acknowledgement}
The present research has been partially sponsored by project $DO-02-275,$
$08.12.2008$ with the Bulgarian NSF.
\end{acknowledgement}


\begin{thebibliography}{99}                                                                                               %


\bibitem {atakhodzaev}M. Atakhodzaev, Ill-posed internal boundary value
problems for the biharmonic equation, VSP, $2002.$

\bibitem {bersjohnschechter}L. Bers, F. John and M. Schechter, Partial
Differential Equations, Interscience, New York, 1964.

\bibitem {browder}Felix Browder, Approximation by Solutions of Partial
Differential Equations, American Journal of Mathematics, vol. 84, no. 1, p.
134, 1962; Functional analysis and partial differential equations. II,
Mathematische Annalen, vol. 145, no. 2, pp. 81-226, 1962. 

\bibitem {hoermander1}Hormander, L., The Analysis of Linear Partial
Differential Operators I, Springer-Verlag, $1990.$ 

\bibitem {kreinnudelman}Krein, M. , Nudelman A., The Markov Moment problem and
Extremal problems, AMS translation from the Russian edition of $1973.$

\bibitem {karlinstudden}S. Karlin, W. Studden, Tchebycheff Systems: with
applications in analysis and statistics, Intersci. Publ., $1968.$

\bibitem {okbook}Kounchev, O., Multivariate Polysplines, Academic Press, San
Diego, $2001.$ 

\bibitem {lax}P. Lax, A stability theorem for solutions of abstract
differential equations and its application to the study of local behavior of
solutions of elliptic equations, Comm. Pure Appl. Math. 9 (1956) 747--766.

\bibitem {lionsmagenes}J.-L. Lions and E. Magenes, Non-homogeneour Boundary
Value Problems and Applications, Springer, Berlin-Heidelberg, $1970.$

\bibitem {malgrange}B. Malgrange, Existence et approximation des solutions des
equations aux derivees partielles et des equations de convolution, Ann. Inst.
Fourier (Grenoble) 6 (1955--1956) 271--355.

\bibitem {schumaker}L. Schumaker, \emph{Spline Functions: basic theory},
Academic Press, NY, $1983.$

\bibitem {tarkhanov}Tarkhanov, N., Approximation on compact sets by solutions
of systems with surjective symbol, Russian Mathematical Surveys , 48(5) (1993)
: 103. 
\end{thebibliography}
\end{document}